\documentclass[xypic,amscd,syntonly,amssymb]{amsart}

\usepackage{graphicx}
\usepackage[all]{xy}
\usepackage{amssymb}
\usepackage{amsmath}

\usepackage{epsfig}
\theoremstyle{plain}
\newtheorem{Thm}{Theorem}
\newtheorem{Prop}[Thm]{Proposition}
\newtheorem{Cor}[Thm]{Corollary}

 \theoremstyle{definition}
\theoremstyle{remark}

\numberwithin{equation}{section}

\begin{document}
 \title{Critical Values of Moment Maps
   on Quantizable Manifolds}

 \author{ ANDR\'{E}S   VI\~{N}A}
\address{Departamento de F\'{i}sica. Universidad de Oviedo.   Avda Calvo
 Sotelo.     33007 Oviedo. Spain. }
 \email{vina@uniovi.es}
\thanks{}

 \maketitle
\begin{abstract}

Let $M$ be a  quantizable symplectic manifold acted on by
$T=(S^1)^r$
  in a Hamiltonian fashion and $J$ a moment map for this action.
  Suppose that
  the set $M^{T}$ of fixed points is discrete and denote by ${\alpha}_{pj}\in{\mathbb Z}^r$ the weights of the isotropy representation
  at $p$. By means of the  $\alpha_{pj}$'s
we define a partition ${\mathcal Q}_+$, ${\mathcal Q}_-$ of $M^T$.
(When $r=1$,
  ${\mathcal Q}_{\pm}$ will be the set of fixed points
such that the half of the Morse index of $J$ at them is even
  (odd)). We prove the existence of a map $\pi_{\pm}:{\mathcal
  Q}_{\pm}\to{\mathcal Q}_{\mp}$ such that
  $J(q)-J(\pi_{\pm}(q))\in I_{\mp}$, for all $q\in {\mathcal
  Q}_{\pm}$, where $I_{\pm}$ is the lattice generated by  the $\alpha_{pj}$'s with $p\in{\mathcal Q}_{\pm}.$
   We define partition functions $N_p$ similar to the ones of
   Kostant \cite{Gui} and we prove that $\sum_{p\in{\mathcal
   Q}_+}N_p(l)=\sum_{p\in{\mathcal
   Q}_-}N_p(l)$, for any $l\in{\mathbb Z}^r$ with $|l|$
   sufficiently large.

\end{abstract}
   \smallskip



\section {Introduction} \label{S:intro}

Let $(M,\,\omega)$    be a closed oriented symplectic manifold of
dimension $2n$. Henceforth we assume that $\frac{1}{2\pi}[\omega]$
belongs to the image of the map $H^2(M,\,{\mathbb Z})\to
H^2(M,\,{\mathbb R})$; in other words, $(M,\,\omega)$ is {\em
quantizable} \cite{W}.
  Let $(L,\,\nabla)$ be prequantum data on $M$, that
is, $L$ is a Hermitian line bundle on $M$ and $\nabla$ is a
connection on $L$ whose curvature is $-i\omega$.

Let $T$ be the torus $(U(1))^r$ and we assume that $T$ acts on $M$
in a Hamiltonian fashion. We denote by $\mu$
$$\mu:M\to {\frak t}^*= \overbrace{i{\mathbb R}\oplus\dots\oplus i{\mathbb R}}^r  $$
the corresponding normalized moment map (that is, we suppose  that
$\int_M\langle\mu,\,X\rangle\omega^n=0$, for all $X\in{\frak t}$).
We will write  $J$ for the ${\mathbb R}^r$ valued map $-i\mu$.

Throughout this Section we assume that the {\em prequantum data
are
 $T$-invariant} (see \cite{G-S} and Section \ref{SectEquivariant}), and that the set of fixed
points  {\em $M^T$ is a set of isolated points.} If $p\in M^T$, we
denote by $\alpha_{pj}\in{\mathbb Z}^r$, $j=1,\dots, n$ the
weights of the isotropy representation $R$ of $T$ on the tangent
space to $M$ at $p$.

 Now we restrict ourselves to
the case $r=1$. Given $p\in M^{U(1)}$, we set
 $${\mathcal A}_p:=\{i \,|\, \alpha_{pi}>0 \},\;\; {\mathcal B}_p:=\{k \,|\, \alpha_{pk}<0
 \},\; \; \sigma(p):=(-1)^{\sharp {\mathcal B}_p }.$$

We put $b:=\sharp{\mathcal B}_p$ and $a:=n-b$. Given a natural
number $l$, we denote by $N_p(l)$ the following cardinal
\begin{equation}\label{Np(l)}
\sharp\Big\{(m_1,\dots,m_a,n_1,\dots,n_b)\,\big|\,
J(p)+\sum_{i\in{\mathcal A}_p}m_i\alpha_{pi}-\sum_{k\in {\mathcal
B}_p} n_k\alpha_{pk}=l,\; m_i\in{\mathbb N}_{>0},\; n_k\in{\mathbb
N}\Big\}.
\end{equation}
  That is, $N_p(l)$ is the number
of times we can write $l$ as a sum of the type specified in
(\ref{Np(l)}). If $a=0$, $N_p$ is the Kostant partition function
associated with the representation $R^*$, but in general our
``partition" functions do not agree with the ones of Kostant.

On the other hand, one has the partition $M^{U(1)} ={\mathcal
Q}_+\cup {\mathcal Q}_{-}$, where
$${\mathcal Q}_{\pm} =\{p\in
M^{U(1)}\,|\,\sigma(p)={\pm 1}\}.$$
In Subsection
\ref{SectionU(1)} we will prove the following Theorem
\begin{Thm}\label{ThmGeneral} For $l\in{\mathbb N}$ sufficiently
large
$$\sum_{p\in{\mathcal Q}_+}N_p(l)=\sum_{p\in{\mathcal
Q}_-}N_p(l).$$
\end{Thm}

Let us denote by  $I_{\pm}$ the ideal of ${\mathbb Z}$ generated
by the set
   $$\{\alpha_{qj}\;|\; j=1,\dots,n ;\;   q\in {\mathcal Q}_{\pm}\}.$$
In Subsection \ref{SectionU(1)} the following Theorem is proved

\begin{Thm}\label{ThmU(1)}
Given $p_{\pm}\in{\mathcal Q}_{\pm}$, there is $q_{\mp}\in
{\mathcal Q}_{\mp}$ such that
 \begin{equation}\label{Jp+-Jp-}
 J(p_{\pm})-J(q_{\mp})\in I_{\mp}.
 \end{equation}
There is a natural number $c_{\pm}$, $1\leq
c_{\pm}\leq\sharp{\mathcal Q}_{\mp}$ such that
\begin{equation}\label{cpmIpm}
c_{\pm}I_{\pm}\subset I_{\mp}.
\end{equation}
\end{Thm}

In a neighborhood of the fixed point $p$ for the circle action
there are suitable coordinates $x_1,\dots, x_n,y_1,\dots,y_n$,
such that the moment map $J$ in this neighbourhood  is given by
\begin{equation}\label{JMorse}
J(p)+\frac{1}{2}\sum_{j=1}^n\alpha_{pj}(x^2_j+y^2_j).
 \end{equation}
  Thus Theorem \ref{ThmU(1)} gives a relation between the values of
the moment map $J$ at the fixed points and the Hessian of $J$,
relative to those coordinates, at these points.

By (\ref{JMorse}) the index of the critical point $p$ of the Morse
function $J$ is $2(\sharp{\mathcal B}_p)$. We can classify the
points $p\in M^{U(1)}$ in even and odd according to the parity of
$(1/2){\rm index }(p)$. Thus $p$ is even iff $p\in {\mathcal
Q}_+$. Hence we have the following Corollary, in whose statement
we put $I_{{\rm odd}}:=I_-.$
 \begin{Cor}\label{Coreven} For each even critical point $p$ of $J$ there is an
 odd critical point $q$ such that
 $$J(p)-J(q)\in I_{{\rm odd}}.$$
 \end{Cor}

Obviously, the statement obtained from Corollary \ref{Coreven}
exchanging ``even" for ``odd" is also true.

\smallskip

In this note we will prove two types of properties for $(U(1))^r$-
Hamiltonian actions:

{\bf (A)} Properties with a combinatorial content (as the
statement of Theorem \ref{ThmGeneral}), in which are involved the
weights of the isotropy representations.

{\bf (B)} Relations satisfied by the difference of  critical
values of the moment map, as the result stated in Theorem
\ref{ThmU(1)}.

\smallskip

Convexity properties of a moment map $\mu$ for a Hamiltonian
action of a torus $T$ on a symplectic manifold were studied by
Atiyah in \cite{A} and by Guillemin and Sternberg in \cite{G-S-0}.
The basic result they proved is that ${\rm im}\,\mu$ is a convex
polytope in ${\frak t}^*$; it is, in fact, the convex hull of the
points $\mu(p)$, where $p$ is a fixed point of the $T$-action.
This result was generalized by Kirwan in \cite{K} to Hamiltonian
actions of compact connected Lie groups.  Furthermore, Theorem 3
of \cite{G-S-0} describes the ``shape of the vertex" $\mu(p)$ of
the polytope ${\rm im}\,\mu$ in terms of the weights of the
isotropy representation of $T$ on the tangent space to $M$ at $p$.

From the convexity properties of the moment map for actions on
quantizable manifolds one deduces:

1.)   The difference between two critical values of $J$
 belongs to integer
lattice of
${\mathbb R}^r$ (see \cite{G-S-0} and \cite{We}). We will obtain
this result again in Section \ref{SectEquivariant} when we explain
briefly the concept of invariant prequantum data (see
(\ref{Jinteger})).

2.) If $J(p)$ and $J(q)$ are vertices of ${\rm im}\,J$ {\em which
determine and edge} of ${\rm im}\,J$, the local convexity theorem
of \cite{G-S-0} implies that $J(p)-J(q)$ belongs to the lattice
generated by $\{\alpha_{qj}\}_j$.

In this paper not all the fixed points play identical role. They
will be classified
in two classes, say ${\mathcal Q}_+$, ${\mathcal Q}_-$, using
 the weights of the respective isotropy
representations. The properties of type ${\bf A}$ give the
equality between the sum of our partitions functions for the
points of ${\mathcal Q}_+$  and the corresponding one for the
points of ${\mathcal Q}_-$. The proof we will give is elementary
and independent of the known formulas for the
  Kostant partition functions \cite{Brion-V}, \cite{Gui1}.


Our properties of type ${\bf B}$ give differences between critical
values of $J$ belonging to different class. So Corollary
\ref{Coreven}, Theorem \ref{ThmU(1)} and, in general, the
properties of type ${\bf B}$ give
restrictions to some differences between critical values of the
moment map, which can not be deduced from the above convexity
properties.

The proofs of  our results are  based in the equivariant index
formula for the Dirac operator \cite{L-M}, \cite{B-G-V}. Although
this formula has been used in different articles to prove ``that
quantization commutes with reduction"
(see for example \cite{D-G-M-W}, \cite{G}, \cite{J-K},
\cite{Mein}, \cite{Verg} or the research report \cite{SJ}),
 in Section \ref{SectEquivariant} we review this
formula and write it in a form which will be convenient for our
applications.

 The   case of circle actions with only isolated fixed points is considered in
Subsection \ref{SectionU(1)}. We will prove Theorem
\ref{ThmGeneral} and Proposition \ref{Propq(p)};  Theorem
\ref{ThmU(1)} is a straightforward consequence of Proposition
\ref{Propq(p)}.

The case when $M$ is acted on by $T=(U(1))^r$ and $M^T$ is a
discrete set is considered in   Subsection \ref{SectionT}. In this
case a way to define a partition of $M^T$ in two subsets is by
means of a polarizing vector (i.e. a vector non-orthogonal to
$\alpha_{pj}$, for all $p$ and $j$). The corresponding properties
of type {\bf A} and type {\bf B} are the statements of Theorem
\ref{ThmGeneralGeneral} and Theorem \ref{Thm(T)}, respectively. We
also define other more general partitions of $M^T$ which will give
rise to Theorem \ref{ThmKostant} and Theorem \ref{Thmnuevo}, whose
statements are properties of type {\bf A} and type {\bf B},
respectively.
 These theorems are the   generalizations of Theorem
\ref{ThmGeneral} and Theorem \ref{ThmU(1)} to the case in which
$r$ is arbitrary. As a corollary we will prove that there is no an
open half space of ${\mathbb R}^r$   containing all the weights
$\alpha_{pj}$. Although the theorems of Subsection
\ref{SectionU(1)}  are particular cases of the corresponding ones
of Subsection \ref{SectionT}, for the sake of clarity of the
exposition we give a direct proof of them.

Finally, in Section \ref{GeneralGeneral} we analyze the case when
$M^T$ is a non necessarily discrete
 set. The statements of Theorem \ref{ThGeneral3} and Corollary \ref{CorThGeneral3} are a properties of
 type {\bf A}. This theorem in turn generalizes Theorem
 \ref{ThmGeneralGeneral}.
  We will apply Theorem \ref{ThGeneral3} to a
 particular case and we will obtain simple combinatorial relations satisfied by the weights of the isotropy
 representations  (see the first two items
 of Proposition \ref{sharpformula}).

  It is not easy to deduce  general simple properties of type  {\bf
 B}   when the fixed points are not isolated. We consider this point in a particular case in
 the Example at the end of the
 paper. The result   is  the third item of Proposition \ref{sharpformula}.


\section{Equivariant Index Theorem}\label{SectEquivariant}

 Given a loop $\{\phi_t  \}_{t\in[0,1]}$ in the Hamiltonian group
 ${\rm Ham}(M,\omega)$ at the identity (see \cite{Mc-S},  \cite{lP01}), let $Y_t$ be the
 time-dependent vector field defined by
 $$\frac{d\,\phi_t}{dt}=Y_t\circ\phi_t.$$
   If $f_t$ is the corresponding normalized Hamiltonian
   ($df_t=-\iota_{Y_t}\omega$), we can pose the problem of
   determining a family $\{\sigma_t\}$ of sections of $L$
   satisfying
   $$\frac{d\sigma_t}{dt}=-\nabla_{Y_t}\sigma_t-if_t\sigma_t,\;\;\;
   \sigma_0=\sigma,$$
   where $\sigma\in\Gamma(M,\,L)$ is a given section. It turns out
   that $\sigma_1=\kappa(\phi)\sigma$, where the constant
   $\kappa(\phi)$ is given by
   $$\kappa(\phi)={\rm exp}\Big( i\int_S\omega-i\int_0^1f_t(p)dt
   \Big),$$
   $p$ being an arbitrary point of $M$ and $S$ an arbitrary
   $2$-cycle in $M$ whose boundary is the curve
   $\{\phi_t(p)\}_t$ (see \cite{V}). If $p$ is a fixed point, i.e.
   $\phi_t(p)=p$ for all $t$, then
\begin{equation}\label{kappafixed}
 \kappa(\phi)={\rm exp}\Big(  -i\int_0^1f_t(p)dt
   \Big).
   \end{equation}

\smallskip

Suppose that the torus $T=(U(1))^r$ acts on $M$ in a Hamiltonian
fashion.
As we said,  by $\mu$ will be denoted the corresponding normalized
moment map. Given $X\in{\frak t}$, we adopt the following
convention
$$d\langle \mu,\, X\rangle=\omega(X_M,\, .),$$
where $X_M$ is the vector field on $M$ defined by
$$X_M(p)=\frac{d}{d\epsilon}\big(e^{-\epsilon X}\cdot p\big)_{\big |\epsilon=0}.$$

For each $X\in{\frak t}$ one can consider the operator
$${\mathcal P}_X:=\nabla_{X_M}+i\langle\mu,\, X
\rangle:\Gamma(M,\,L)\to \Gamma(M,\,L).$$
 ${\mathcal P}$ is a representation of the Lie algebra ${\frak
 t}$. The prequantum data are said to be $T$-invariant if the representation
 ${\mathcal P}$ can be lifted to  a representation of $T$
 \cite{G-S}.

From now on we assume that {\em the prequantum data
    $(L,\,\nabla)$  are $T$-invariant}. If
    $X$ belongs to the integer lattice $\Lambda$ of ${\frak t}\,$
   ($\Lambda={\rm ker(exp:}\,{\frak t}\to T)$) and  $p$ is a fixed
   point for the $T$-action, by applying (\ref{kappafixed}) to
   $\phi_t(q)=e^{tX}\cdot q$ one obtains
   \begin{equation}\label{Jinteger}
   \langle \mu(p),\,X\rangle\in 2\pi{\mathbb Z}.
   \end{equation}

   \smallskip

   Given an almost complex structure ${\mathcal I}$ on $M$ compatible with
   $\omega$, we can consider the spaces $\Omega^{ij}(M,\,L)$ of
   $L$-valued forms of type $(i,j)$. By means of the Hermitian connection on the canonical bundle on $M$
   we can define   the spin-${\mathbb
   C}$ Dirac operator (see \cite{L-M})
$$ D_{\mathbb C}:\Omega^{0,{\rm even}}(M,\,L)\to \Omega^{0,{\rm
odd}}(M,\,L),$$ and the corresponding virtual vector space
$${\rm ind}(D_{\mathbb C}):=[{\rm ker}(D_{\mathbb C})]-[{\rm coker}(D_{\mathbb
C})].$$

 Since two almost complex structures on $M$ compatible with
 $\omega$ are homotopic, ${\rm ind}(D_{\mathbb C})$ is independent
 of ${\mathcal I}$. If ${\mathcal I}$ is $T$-invariant, there is a natural
 representation of $T$ in ${\rm ind}(D_{\mathbb C})$.

 We will study the character $\chi$ of this
 representation of $T$. If $t\in T$, $\chi(t)$ is a Laurent
 polynomial, that is, a finite sum of the form
 \begin{equation}\label{Laurentpoly}
\chi(t)=\sum_{{\rm finite \; sum}}E(\rho)t^{\rho},
 \end{equation}
where $\rho\in{\mathbb Z}^r$ and $E(\rho)\in {\mathbb Z}$.

One defines the moment $\mu^L$ as the $T$-equivariant form with
values in ${\rm End}\,(L)$  given by the expression (see
\cite{B-G-V})
$$\mu^L(X)={\mathcal P}_X-\nabla_{X_M},$$
${\mathcal P}_X$ being the above infinitesimal action on
$\Gamma(M,L)$. The equivariant curvature of the connection
$\nabla$ is $R_T(X)=R+i\mu^L(X)$, where $R$ is the curvature of
$\nabla$. So $R_T(X)=-i\omega-\langle\mu,\,X\rangle$ and this is a
$T$-equivariant closed form; i.e.,
$$(d_T\,R_T)(X):=d(R_T(X))+i\iota(X_M)(R_T(X))=0.$$

The $T$-equivariant Chern character of $L$, ${\rm Ch}_T(L,X)$ is
defined as the equivariant cohomology class of
 \begin{equation}\label{Cherncharac}
 {\rm exp}(iR_T(X))=
{\rm exp}(\omega-i\langle\mu,\,X\rangle).
\end{equation}
  Similarly
we can define  ${\rm Td}_T(M,X)$, the $T$-equivariant Todd class
of $M$.

 For $X\in {\frak t}$ sufficiently close to zero so that $M_0$, the zero
 set of $X_M$, is $\{p\in M\,:\,e^X\cdot p=p\}$, the equivariant
 index formula  gives the following expression   for
 $\chi(e^X)$
 \begin{equation}\label{chi(ex)}
 \chi(e^X)=\int_M {\rm Ch}_T(L,X){\rm Td}_T(M,X).
\end{equation}
The localization formula in equivariant cohomology gives (see
\cite{A-S}, \cite{B-G-V})
\begin{equation}\label{chi(ex)2}
\chi(e^X)=\int_{M_0} \frac{{\rm Ch}_T(L,X){\rm
Td}_T(M,X)}{e_T({\mathcal N},X)},
 \end{equation}
where $e_T({\mathcal N},X)$ is the equivariant Euler class of the
normal bundle ${\mathcal N}$ to $M_0$ in $M$ \cite{B-G-V}.

 Now we assume that   the $T$-action has only
isolated points. If $M_0$ is equal to $M^T=\{p\in M\,|\,t\cdot
p=p,\;\text{for all }\; t\in T \}$ it follows from
 (\ref{Cherncharac})  together with (\ref{chi(ex)2}) the following
 expression for the character $\chi$

\begin{equation}\label{chi(ex)3}
\chi(t=e^X)=\sum_{p\in M^T} \frac{t^{J(p)}}
{\prod_{j=1}^n\big(1-t^{-\alpha_{pj}} \big)},
\end{equation}
where $\alpha_{pj}=(\alpha^1_{pj},\dots,\alpha^r_{pj} )\in{\mathbb
Z}^r$, $j=1,\dots,n$ are the weights of the $T$-action on $T_pM$,
and
 $$t^{J(p)}=e^{-i\langle\mu(p),\,X\rangle}.$$


\section{Isolated critical points.}\label{SectionPropert}
 Throughout this Section {\em we will assume that $M^T$ is a discrete
 set}.

\subsection{Case T=U(1)}\label{SectionU(1)}

Now let us suppose that $r=1$. Then for a generic  $z\in U(1)$, as
$M^{U(1)}$ is a discrete set
 \begin{equation} \label{ChiS1}
 \chi(z)=\sum_{p\in
M^{U(1)}}\frac{z^{J(p)}}{\prod_{j=1}^n\big(1-z^{-\alpha_{pj}}
\big)}.
\end{equation}
 By (\ref{Jinteger}) $J(p)\in {\mathbb Z}$, and  the right hand side
 of (\ref{ChiS1})  extends  to a meromorphic function on
 ${\mathbb C}\cup\{\infty\}$.

\smallskip

{\bf Proof of Theorem \ref{ThmGeneral}.}

 If $0<|z|<1$ and $a\in{\mathbb Z}$
\begin{equation}
\frac{1}{1-z^{-a}}=\begin{cases}
-\sum_{m=0}^{\infty}z^{(m+1)a},\;\;\text{if}\,\;a>0  \\

\\

\sum_{m=0}^{\infty}z^{-ma},\;\;\text{if}\,\;a<0
\end{cases}
\end{equation}
Therefore, if $0<|z|<1$ the righthand side of (\ref{ChiS1}) can be
written
 \begin{equation}\label{RightHand}
(-1)^n\sum_{p\in M^{U(1)}}\sigma(p)z^{J(p)}\Big(\prod_{i\in
{\mathcal A}_p}\sum_{m=1}^{\infty}z^{m\alpha_{pi}}   \Big)
\Big(\prod_{k\in {\mathcal B}_p}
\sum_{m=0}^{\infty}z^{-m\alpha_{pk}}   \Big).
 \end{equation}
So an arbitrary summand of (\ref{RightHand}) has the following
form
 \begin{equation}\label{sigmz}
 (-1)^n\sigma(p)z^{{J(p)}+\sum_i m_i\alpha_{pi}-\sum_k
n_k\alpha_{pk}},
 \end{equation} with $n_k\in{\mathbb N}$ and $m_i \in{\mathbb N}_{>0}$ .
The number of times that the monomial $+(-1)^nz^l$ appears in
(\ref{RightHand}) is precisely
 $$\sum_{p\in{\mathcal Q}_{+}}N_p(l),$$
 and similarly for the monomial $-(-1)^nz^l$.
 Since $\chi(z)$ is  a {\em finite} sum of
 monomials in the variable $z$, all except a finite number of terms of the form
 (\ref{sigmz}) with $\sigma(p)=+1$,
must cancel  all the terms with $\sigma(p)=-1$ except a finite
number of them. That is, for $l\in{\mathbb N}$ large enough
$\sum_{p\in{\mathcal Q}_+} N_p(l)=\sum_{p\in{\mathcal Q}_-}
N_p(l).$

\qed

\smallskip

 {\it Example 1 ($S^1$-action on ${\mathbb C}P^2$).}

 Now we consider the symplectic toric manifold $M$ asssociated to
 the Delzant polytope $\Delta$ in $({\mathbb R}^2)^*$ with vertices
 $P(0,0)$, $Q(1,0)$ and $R(0,1)$. By $J_T$ we denote the moment
 map of the ${\mathbb T}^2$-action. Let us consider the
 $S^1$-action generated by $X=(x,y)\in{\mathbb Z}^2$ through the
 ${\mathbb T}^2$-action. The corresponding moment map is
 $J=\langle J_T,\, X\rangle$. If $x,\, y>0$ and $x>y$, then
 $M^{S^1}=\{p,q,r\}$, where $p=J_T^{-1}(P)$,\;
 $q=J_T^{-1}(Q)$,\;$r=J_T^{-1}(R)$; moreover $J(p)=0$,\;$J(q)=x$ and $J(r)=y$.
 The weights of the isotropy
 representations of $U(1)$ are
 $$\alpha_{p1}=x,\;
 \alpha_{p2}=y;\,\;\alpha_{q1}=-x,\;\alpha_{q2}=-x+y;\,\;
 \alpha_{r1}=x-y,\;\alpha_{r2}=-y.$$
 Thus ${\mathcal A}_p=\{1,2\},\,{\mathcal A}_q=\emptyset,\,{\mathcal
 A}_r=\{1\}$, $\,{\mathcal Q}_+=\{p,q\}$ and ${\mathcal Q}_-=\{r\}$.

 Given $l\in{\mathbb N}$ sufficiently large, we set
 $${\mathcal S}_p:=\{(m_1,\,m_2)\,|\, l=m_1x+m_2y  \},\;\;
{\mathcal S}_q:=\{(n_1,\,n_2)\,|\, l=x+n_1x+n_2(x-y)  \}$$
$${\mathcal S}_r:=\{(m,\,n)\,|\, l=y+m(x-y)+ny  \},$$
where $m,m_1,m_2\in{\mathbb N}_{>0}$ and $n,n_1,n_2\in{\mathbb
N}$. The map ${\mathcal S}_r\to {\mathcal S}_p\coprod{\mathcal
S}_q$ defined by
$$
(m,n)\mapsto\begin{cases} (m,\,1-m+n)\in{\mathcal S}_p,\;\,\text{if}\,\, 1-m+n>0 \\
                           (n,\,-1+m-n)\in{\mathcal S}_q,\;\,\text{if}\,\, 1-m+n\leq
                           0.
                           \end{cases}
 $$
is bijective. As $N_b(l)=\sharp{\mathcal S}_b$, for $b\in\{p,q,r
\}$, it follows $N_r(l)=N_p(l)+N_q(l)$, which is the assertion of
Theorem \ref{ThmGeneral} in this particular case.

\smallskip

\begin{Prop}\label{Propq(p)}
Given $p\in{\mathcal Q}_+$, there is $q(p)\in {\mathcal Q}_-$ such
that
 $$J(p)-J(q(p))\in I_-.$$
   Given
$i\in{\mathcal A}_p$, there is a natural number $c$, $1\leq
c\leq\sharp{\mathcal Q}_-$ such that $c\alpha_{pi}\in I_-$.
\end{Prop}

{\bf Proof.}

  As (\ref{RightHand}) is a finite sum,  given $p\in {\mathcal Q}_+$, $\tilde i\in {\mathcal A}_p$ and
 $m\in{\mathbb N}$ sufficiently large  there exist
$q\in {\mathcal Q}_-$, $m'_i,n'_k$ such that
$$J(p)+m\alpha_{p\tilde i}=J(q)+\sum_i m'_i\alpha_{qi}-\sum_k
n'_k\alpha_{qk}.$$
  So
  \begin{equation}\label{Jp-JqI}
   J(p)-J(q)+m \alpha_{p\tilde i}\in I_-.
  \end{equation}

   If we repeat the argument for $m+1$ we will conclude that there
   is $q_1\in {\mathcal Q}_{-}$ such that
$J(p)-J(q_1)+(m +1)\alpha_{p\tilde i}\in I_-$, etc. Since the
number of points of ${\mathcal Q}_{-}$ is finite, there is a point
in ${\mathcal Q}_{-}$, which we denote by $\Hat q$, and an integer
$c$, $1\leq c \leq \sharp {\mathcal Q}_{-}$ that
 $$J(p)-J(\Hat q)+m \alpha_{p\tilde i}\in I_-\; \text{and}\; J(p)-J(\Hat q)+(m +c)\alpha_{p\tilde i}\in
 I_-,$$
 for some $m$ large enough. Thus
 \begin{equation}\label{calpphaI-}
  c\alpha_{p\tilde i}\in I_-.
 \end{equation}

 Reasoning  with the element  $c\alpha_{p\tilde i}$ as before with $\alpha_{p\tilde
 i}$, we can conclude that there is an element $q(p)\in {\mathcal Q}_-$ such that
$$J(p)-J(q(p))+ m_1 c\alpha_{p\tilde i}\in I_-,$$
with $m_1$ a natural number sufficiently large. From
(\ref{calpphaI-})
 it follows that
 $$J(p)-J(q(p))\in I_-.$$
\qed
\smallskip

{\bf Proof of Theorem \ref{ThmU(1)}.}
 The
arguments used in the proof of Proposition \ref{Propq(p)} applied
to the case that $p\in{\mathcal Q}_-$, together with Proposition
\ref{Propq(p)} complete the proof of Theorem \ref{ThmU(1)}.
  \qed

\smallskip


\subsection{Case $T=(U(1))^r$}\label{SectionT}

Let $t=(t_1,\dots,t_r)$ be a generic element of $T$ sufficiently
close to the identity. As we are assuming that $T$-action has only
isolated points, then $\chi(t)$ can be calculated by
(\ref{chi(ex)3}). The right hand side of (\ref{chi(ex)3}) is
defined on the region  of ${\Bbb C}^r$ consisting of the points
$t$ such that $t^{\alpha_{pj}}\ne 1$ for  any $p\in M^T$ and any
$j$.

As $\alpha_{pj}\ne 0$, for a generic vector $u\in{\mathbb R}^r$
$$\alpha_{pj}(u):=\sum_{e=1}^ru_e\alpha^e_{pj}\ne 0,$$
for all $p\in M^T$ and all $j=1,\dots,n$. Such a vector is said to
be polarizing. By means of $u$ we will define a partition of $M^T$
of two subsets.

We denote
 $$ {\mathcal A}_p(u):=\{i \,|\, \alpha_{pi}(u)>0 \},\;\; {\mathcal B}_p(u):=\{k \,|\,
 \alpha_{pk}(u)<0\},\;\;
 J_p(u)=\sum_e u_eJ^e_p. $$
Given a  number $l\in{\mathbb N}u_1+\dots +{\mathbb N}u_r$, one
defines $N_p(l,u)$ by the formula obtained from (\ref{Np(l)})
exchanging $J(p)$ for $J_p(u)$, $\alpha_{pj}$ for
$\alpha_{pj}(u)$, ${\mathcal A}_p$ for ${\mathcal A}_p(u)$ and
${\mathcal B}_p$ for ${\mathcal B}_p(u)$.

Similarly we define
 $$\sigma_{p}(u):=(-1)^{\sharp {\mathcal B}_p(u) }, \; \; {\mathcal Q}_{\pm}(u)=\{p\in M^T,|\,\sigma_{p}(u)=\pm 1  \}. $$

Fixing a branch of the logarithmic function, then
  $v(\lambda):=(\lambda^{u_1},\dots,\lambda^{u_r})\in{\mathbb C}^r$ is a single-valued analytic
   function on a region ${\mathcal R}$
   of ${\mathbb C}\setminus\{0\} $.  Let $\lambda$ be a point of  ${\mathcal R}$
   with
   $0<|\lambda|<1$, then
   $$v^{-\alpha_{pj}(u)}:=\prod_e(\lambda^{u_e})^{-\alpha^e_{pj}}=\lambda^{-\alpha_{pj}(u)}.$$
The right hand side of (\ref{chi(ex)3}) for $t=v$ is
\begin{equation}\label{RightHandGeneral}
(-1)^n\sum_{p\in M^{T}}\sigma_p(u)\lambda^{J_p(u)}\Big(\prod_{i\in
{\mathcal A}_p(u)}\sum_{m=1}^{\infty}\lambda^{m\alpha_{pi}(u)}
\Big) \Big(\prod_{k\in {\mathcal B}_p(u)}
\sum_{m=0}^{\infty}\lambda^{-m\alpha_{pk}(u)}   \Big),
 \end{equation}
By (\ref{Laurentpoly}) the coefficient of $\lambda^l$ must be
zero,  for $l\in{\mathbb N}u_1+\dots +{\mathbb N}u_r$ and $|l|$
large enough. As in the proof of Theorem \ref{ThmGeneral}, this
coefficient
$$(-1)^n\big(\sum_{p\in{\mathcal Q}_+(u)}N_p(l,u)-\sum_{p\in{\mathcal
Q}_-(u)}N_p(l,u)\big).$$
 Thus we have proved the
following Theorem
 \begin{Thm}\label{ThmGeneralGeneral}
  Let $u$ be a vector of the Euclidean space ${\mathbb R}^r$ which is
  non-orthogonal to
$\alpha_{pj}$ for all $p\in M^T$ and all $j\in\{1,\dots, n\}$, and
$l\in{\mathbb N}u_1+\dots +{\mathbb N}u_r$ with $|l|$ sufficiently
large, then
$$\sum_{p\in{\mathcal Q}_+(u)}N_p(l,u)=\sum_{p\in{\mathcal
Q}_-(u)}N_p(l,u).$$

 \end{Thm}

From the proof of Theorem \ref{ThmGeneralGeneral}  we deduce the
following Corollary

\begin{Cor}\label{Corohalfspace}
There is no an open half space in ${\mathbb R}^r$ which contains
all the $\alpha_{pj}$.
\end{Cor}
{\bf Proof.} Suppose the Corollary were false. Let $0\ne u$ be a
vector in ${\mathbb R}^r$ orthogonal to the hyperplane boundary of
the half space. Then whether ${\mathcal A}_p(u)$ or ${\mathcal
B}_p(u)$ would be empty for all $p$. Thus all the summands in
(\ref{RightHandGeneral}) would have the same sign. This
contradicts the fact that  (\ref{RightHandGeneral}) is  a finite
sum.
 \qed

\smallskip

As in the proof of Proposition \ref{Propq(p)},
  given $p\in{\mathcal Q}_+(u)$,
$i\in{\mathcal A}_p(u)$ and $m>>1$, there exists $q_u\in{\mathcal
Q}_-(u)$ such that
$$ J_p(u)+m\alpha_{pi}(u)=J_{q_u}(u)+\sum_im'_i\alpha_{q_u i}(u)-\sum_kn'_k\alpha_{q_u k}(u). $$
On the other hand, for $u'\in{\mathbb R}^r$ sufficiently close to
$u$ we have
$$ {\mathcal A}_p(u)={\mathcal A}_p(u') ,\;\; {\mathcal B}_p(u)=\;\; {\mathcal B}_p(u')
, \; \; {\mathcal Q}_{\pm}(u)= \; \; {\mathcal Q}_{\pm}(u').
$$
Since ${\mathcal Q}_-(u)$ is a finite set, for $u'$ in a
neighborhood of $u$, the corresponding $q_{u'}\in {\mathcal
Q}_-(u')={\mathcal Q}_-(u)$ can be taken equal to $q_u$. Therefore
we have for  any $u'$ in a neighborhood of $u$
$$ J_p(u')+m\alpha_{pi}(u')=J_{q_u}(u')+\sum_im'_i\alpha_{q_u i}(u')-\sum_kn'_k\alpha_{q_u k}(u'). $$
That is,
$$ \sum_e u'_e\big(J_p^e+m\alpha^e_{pi}-J^e_{q_u}-\sum_im'_i\alpha^e_{q_u i}+\sum_kn'_k\alpha^e_{q_u k}\Big)=0. $$
Thus for each $e=1,\dots, r$
 \begin{equation}\label{Ja(p)}
 J^e_p-J^e_{q_u}+m\alpha^e_{pi}\in I^e_-(u),
 \end{equation}
$I^e_{\pm}(u)$ being the ideal of ${\mathbb Z}$ generated by
$\{\alpha^e_{qj}\,|\, q\in{\mathcal Q}_{\pm}(u),\, j=1,\dots,n
\}$.

We have obtained a formula similar to (\ref{Jp-JqI}). Reasoning as
in Subsection \ref{SectionU(1)} one can prove the following
Theorem
 \begin{Thm}\label{Thm(T)} Let $u$ be a polarizing vector of the Euclidean space ${\mathbb
 R}^r$.
Given $p_{\pm}\in{\mathcal
 Q}_{\pm}(u)$ there exists $q_{\mp}\in {\mathcal Q}_{\mp}(u)$ such
 that
 $$J^e(p_{\pm})-J^e(q_{\mp})\in I_{\mp}^e(u),$$
 for all $e\in\{1,\dots, r\}$.

 Moreover, there is a natural number $c^e_{\pm}$, with $1\leq c^e_{\pm}\leq
 {\mathcal Q}_{\mp}(u)$, such that
 $$c^e_{\pm}I^e_{\pm}(u)\subset I^e_{\mp}(u).$$
 \end{Thm}

\smallskip

For $p\in M^T$ and $j\in\{1,\dots,n\}$ we write
$${\mathcal R}_{pj}^{\pm}=\big\{ (x_1,\dots,x_r)\in{\mathbb R}^r \,|\,\pm\sum_ex_e\alpha^e_{pj}>0  \big\}.$$
Let $t=(t_1,\dots,t_r)\in{\mathbb C}^r$, such that
$(\log|t_1|,\dots,\log|t_r|)\in {\mathcal R}^+_{pj},$ then
$|t^{-\alpha_{pj}}|<1$. Hence
$$\frac{1}{1-t^{-\alpha_{pj}}}=\sum_{m\geq 0}t^{-m\alpha_{pj}}.$$
If $(\log|t_1|,\dots,\log|t_r|)\in {\mathcal R}^-_{pj},$ then
$$\frac{1}{1-t^{-\alpha_{pj}}}=-\sum_{m\geq 1}t^{m\alpha_{pj}}.$$

We denote by  $\epsilon$ a map $M^T\times \{1,\dots, n\}\to\{\pm
1\}$, such that
 \begin{equation}\label{Bigintersec}
{\mathcal
K}(\epsilon):=\bigcap_{p,j}\mathcal{R}^{\epsilon(p,j)}_{pj}\ne\emptyset
.
 \end{equation}
 Since each $\mathcal{R}^{\epsilon(p,j)}_{pj}$ is an open half space of ${\mathbb R}^r$,
if (\ref{Bigintersec}) holds,  this intersection contains
infinitely many points of ${\mathbb R}^r$ as close to $0$ as we
wish.

 We put
$${\mathcal A}_p(\epsilon)=\{j\,|\,\epsilon (p,j)=-1\},\;\;
{\mathcal B}_p(\epsilon)=\{j\,|\,\epsilon (p,j)=1\},\;\;
\sigma(p,\epsilon)=(-1)^{\sharp{\mathcal A}_p(\epsilon)}.$$

 For any
$t=(t_1,\dots,t_r)\in {\mathbb C}^r$ close to the identity of $T$
and such that
 $$(\log|t_1|,\dots,\log|t_r|)\in {\mathcal K}(\epsilon)$$
  the right hand side of (\ref{chi(ex)3}) is
\begin{equation}\label{RightHandepsilon}
 \sum_{p\in
M^{U(1)}}\sigma(p,\epsilon)t^{J(p)}\Big(\prod_{i\in {\mathcal
A}_p(\epsilon)}\sum_{m\geq 1} t^{m\alpha_{pi}}   \Big)
\Big(\prod_{k\in {\mathcal B}_p(\epsilon)} \sum_{m\geq0}
t^{-m\alpha_{pk}}   \Big).
 \end{equation}

  Given $l\in{\mathbb Z}^r$
we denote by $N_p(l,\epsilon)$
$$\sharp\Big\{(m_1,\dots,m_a,n_1,\dots,n_b)\,\big|\,
J(p)+\sum_{i\in{\mathcal A}_p(\epsilon)}m_i\alpha_{pi}-\sum_{k\in
{\mathcal B}_p(\epsilon)} n_k\alpha_{pk}=l,\; m_i\in{\mathbb
N}_{>0},\; n_k\in{\mathbb N}\Big\}.$$
 Where $a=\sharp {\mathcal A}_p(\epsilon)$ and $b=\sharp {\mathcal
 B}_p(\epsilon)$.
If $\sharp {\mathcal A}_p(\epsilon)=0$, then $N_p(l,\epsilon)$ is
 the Kostant partition function corresponding to the representation
 of $T$ with weights $-\alpha_{pk}$ (see \cite{Gui}).

Now we define the partition  ${\mathcal Q}_{+}(\epsilon)$,
${\mathcal Q}_{-}(\epsilon)$ of $M^T$, where
 $${\mathcal
Q}_{\pm}(\epsilon)=\{p\in M^T\,|\,\sigma(p,\epsilon)=\pm 1\}.$$ As
in  the preceding Subsection, given $l\in{\mathbb Z}^r$ the
 number of times that the monomial $+t^l$ (resp. $(-1)t^l$) appears in
 (\ref{RightHandepsilon}) is
 $$\sum_{p\in{\mathcal Q}_+(\epsilon)}N_p(l,\epsilon)\;\;\;
 \text{(resp.}\; \sum_{p\in{\mathcal Q}_-(\epsilon)}N_p(l,\epsilon)\text{)} .$$
 Hence for all $l\in{\mathbb Z}^r$ with $|l|$ sufficiently large
$$\sum_{p\in{\mathcal Q}_+(\epsilon)}N_p(l,\epsilon)=
\sum_{p\in{\mathcal Q}_-(\epsilon)}N_p(l,\epsilon) .$$
 We have the following theorem
\begin{Thm}\label{ThmKostant} For any map
$$\epsilon:M^T\times \{1,\dots,n\}\to\{\pm 1\},$$
such that ${\mathcal K}(\epsilon)\ne\emptyset$ and for any
$l\in{\mathbb Z}^r$ with $|l|$ big enough,
$$\sum_{p\in{\mathcal Q}_+(\epsilon)}N_p(l,\epsilon)=
\sum_{p\in{\mathcal Q}_-(\epsilon)}N_p(l,\epsilon) .$$
\end{Thm}

If in the proof of Proposition \ref{Propq(p)} we substitute
${\mathcal A}_p$ by ${\mathcal B}_p(\epsilon)$, ${\mathcal Q}_+$
by ${\mathcal Q}_+(\epsilon)$ and $I_-$ by the  lattice in
${\mathbb R}^r$
$$I_-(\epsilon):=\sum_{p\in{\mathcal Q}_-
(\epsilon)}\sum_j{\mathbb Z}\alpha_{pj},$$
 we obtain a proof of the following Theorem
 \begin{Thm}\label{Thmnuevo} Under the hypotheses of Theorem
 \ref{ThmKostant}, given
 $p\in{\mathcal Q}_+(\epsilon)$, there exists $q\in{\mathcal
 Q}_-(\epsilon)$, such that $J(p)-J(q)\in I_-(\epsilon)$.

 Given $i\in{\mathcal B}_p(\epsilon)$, there exists a natural
 number $c$, $1\leq c\leq\sharp {\mathcal Q}_-(\epsilon)$, such that
  $c\alpha_{pi}\in I_-(\epsilon)$.
  \end{Thm}


\section{General case} \label{GeneralGeneral}

 In the case when $M^{T}$ is not a discrete set  the relations which
 satisfy the critical values of $J$ and  the weights of the
 isotropy representations are more complicate. If $F$ is a connected component of $M^T$, we
 denote by ${\mathcal N}_{F}$ the normal bundle   to $F$ in
 $M$. Let
 $$\oplus_{j=1}^s{\mathcal N}_{Fj}$$
   be a decomposition  of ${\mathcal
 N}_{F}$ as direct sum of on line bundles, such that $T$ acts on
 ${\mathcal N}_{Fj}$ with weight $\alpha_{Fj}\in{\mathbb Z}$, with
 $j=1,\dots, s$ (assumed that $s=(1/2)\,{\rm codim}\, F$).
The value of $\chi(t)$ is given
 by the following expression
(see for example \cite{J-K} or \cite{SJ})
\begin{equation}\label{Chinodiscret0}
\chi(t)=\sum_F \chi_F(t)
\end{equation}
where $F$ runs on the set of components of $M^T$ and
\begin{equation}\label{Chinodiscret1}
\chi_F(t)=t^{J(F)}\int_F \frac{e^{\omega}\,{\rm Td}(F)}
 {\prod_{j=1}^s\big( 1-t^{-\alpha_{Fj}}e^{-c_1( {\mathcal N}_{Fj})     }       \big)  }
\end{equation}

Here we fix a polarizing vector $u\in{\mathbb R}^r$ and we put
$v(\lambda)=(\lambda^{u_1},\dots,\lambda^{u_r})\in{\mathbb C}^r$,
as in the     Subsection \ref{SectionT}, where $0<|\lambda|<1$ is
a point of the region ${\mathcal R}$ mentioned in this Subsection.
 We set
$\alpha_{Fj}:=\sum u_e\alpha^e_{Fj}$ and $J(F):=\sum u_eJ^e(F)$.
 The right hand side of (\ref{Chinodiscret1}) for $t=v$ is
  \begin{equation}\label{chiFT}
\lambda^{J(F)}\int_F \frac{e^{\omega}\,{\rm Td}(F) }
 { {\prod_{j=1}^s\big(
1-\lambda^{-\alpha_{Fj}}e^{-c_1( {\mathcal N}_{Fj})     } \big) }
 }.
 \end{equation}
 We will write $\tau_j:=\lambda^{-\alpha_{Fj}}$ and $\gamma_j:= e^{-c_1( {\mathcal
 N}_{Fj})}-1$. Then
  $$ \frac{1}
  { \prod_{j=1}^s\big(
1- \tau_j(1+\gamma_j)    \big) }=\prod_{j=1}^s\Big(\sum_{n\geq
0}\gamma_j^n\frac{\tau_j^n}{(1-\tau_j)^{n+1}}\Big).
  $$
  So the extension of $\chi_F$ to  $t=v$ is
  \begin{equation}\label{chiFT3}
  \chi_F(t)=\lambda^{J(F)}\sum_{n_1,\dots,n_s\geq 0}A_{n_1\dots
  n_s}(F)
  \frac{\tau_1^{n_1}}{(1-\tau_1)^{n_1+1}}\dots
  \frac{\tau_s^{n_s}}{(1-\tau_s)^{n_s+1}},
   \end{equation}
  where
 \begin{equation}\label{Forman}
  A_{n_1\dots n_s}(F):=\int_F e^{\omega}\,{\rm
  Td}(F)\,\gamma_1^{n_1}\dots\gamma_s^{n_s}.
  \end{equation}

  On the other hand
 \begin{equation}\label{bigfraction}
\frac{\tau^m}{(1-\tau)^{m+1}}=\begin{cases} \sum_{l\geq
0}C_-(m,\,l)\tau^l,\;\;\text{if}\;\; |\tau|<1 \\

    \\

(-1)^{m+1}\sum_{l\geq 0}C_+(m,\,l)\tau^{-l},\;\;\text{if}\;\;
|\tau|>1,
\end{cases}
\end{equation}
 where
$$
C_-(m,\,l)=\sharp\{(a_1,\dots,a_{m+1})\,|\, a_j\in{\mathbb N},\;
m+a_1+\dots+a_{m+1}=l \}
$$
$$
C_+(m,l)=\sharp\{(a_1,\dots,a_{m+1})\,|\, a_j\in{\mathbb N},\;
a_1+\dots+a_{m+1}=l-1 \}.
$$
We will put $\sigma_{Fj}={\rm sign}\,(\alpha_{Fj})$, and $\tilde
C_{\sigma_{ Fj}}(m,\,l):=
(-\sigma_{Fj})^{m+1}C_{{\sigma}_{Fj}}(m,\,l),$ then the fraction
(\ref{bigfraction}) when $\tau=\tau_j$ is equal to
$$\sum_{l\geq 0}  \tilde C_{{\sigma}_{Fj}}(m,\,l)\tau_j^{-\sigma_{Fj}l}
$$
 If we define
 $$\tilde{\alpha}_{Fj}=\begin{cases} \alpha_{Fj},\;\;{\rm if}\;
 \alpha_{Fj}>0  \\
 -\alpha_{Fj},\;\;{\rm if}\;
 \alpha_{Fj}<0.
 \end{cases} $$
It follows from (\ref{chiFT3})
\begin{equation}\label{chiFT1}
\chi_F(t)=\lambda^{J(F)}
 \sum_{n_1,\dots,n_s}A_{n_1,\dots,n_s}\sum_{l_1,\dots,l_s\geq
 0}\tilde C_{\sigma_{F1}}(n_1,\,l_1)\dots \tilde
 C_{\sigma_{Fs}}(n_s,\,l_s)\lambda^{\sum_j\tilde \alpha_{Fj}l_j}.
 \end{equation}
 As $A_{n_1,\dots,n_s}(F)$ vanishes if there exists $j$ with
 $n_j>(1/2){\rm dim}\,M$, the first sum in (\ref{chiFT1}) can be
 restricted to $\vec n=(n_1,\dots,n_{s_F})\in{\mathcal G}(F)$,
 where $s_F:=(1/2)\,{\rm codim}\,F$,
 $${\mathcal G}(F)=\{\vec n\in\overbrace{G\times\dots\times G}^{s_F}  \},$$
 and $G:=\{1,\dots,n\}$. The coefficient of $\lambda^k$ in
 $\chi(t)$ is
 $$\sum_F\Big( \sum_{\vec n\in {\mathcal G}(F)}A_{\vec n}(F)\sum_{l\in{\mathcal L}(F,k)}
  \big( \prod_{j=1}^{s_F}\tilde C_{\sigma_{Fj}}(n_j,l_j)  \big) \Big),$$
where
$${\mathcal L}(F, k)=\{(l_1,\dots,l_{s_F})\,|\, l_j\in{\mathbb N},\;\sum_j\tilde
\alpha_{Fj}l_j=k-J(F)  \}.$$

For   $k$ large enough the coefficient of $\lambda^k$ in $\chi(t)$
is zero. Thus we have proved the following Theorem
\begin{Thm}\label{ThGeneral3}
For any polarizing vector $u\in{\mathbb R}^r$ and $k\in {\mathbb
N}u_1+\dots+{\mathbb N}u_r$ with $|k|$ sufficiently large
 \begin{equation}\label{Formula}
\sum_F\Big( \sum_{\vec n\in {\mathcal G}(F)}A_{\vec
n}(F)\sum_{l\in{\mathcal L}(F,k)}
  \big( \prod_{j=1}^{s_F}\tilde C_{\sigma_{Fj}}(n_j,l_j)  \big)
  \Big)=0.
   \end{equation}
\end{Thm}

\smallskip

 We denote by  $D(F,\vec n,k)$ for the combinatorial number
$$ \sum_{ l\in {\mathcal
 L}(F,k)}
 \big(\prod_{j=1}^{s_F}  C_{\sigma_{Fj}}(n_j,\,l_j)\big).$$
  Let
$$\tau(F,\vec n):=\prod_{j=1}^{s_F} (-\sigma_{Fj})^{n_j+1}$$
and
$${\mathcal H}:= \{(F,\vec n)\,|\,F\in{\mathcal F},\;\;\vec
n\in {\mathcal G}(F)\}.$$
 On ${\mathcal H}$ we have    the map
 $\tau:(F,\vec n)\in{\mathcal H}\mapsto \tau(F,\vec
 n)\in\{\pm 1\}$. One has the obvious partition ${\mathcal H}^+$,
 $\, {\mathcal H}^-$ for ${\mathcal H}$. From Theorem
 \ref{ThGeneral3} we deduce the following Corollary whose
 statement has the form of that of Theorem
  \ref{ThmKostant}.
 \begin{Cor}\label{CorThGeneral3}
  For any
 $k\in {\mathbb
N}u_1+\dots+{\mathbb N}u_r$ with $|k|$ sufficiently large and
  any generic vector
  $u\in{\mathbb R}^r$
  $$\sum_{(F,\vec n)\in {\mathcal H}^+}A(F,\vec n)D(F,\vec n,k)=
  \sum_{(F,\vec n)\in {\mathcal H}^-}A(F,\vec n)D(F,\vec n,k),$$
where $A(F,\vec n)$ is given by (\ref{Forman}).
 \end{Cor}

\smallskip

{\it Remark.} If $M^T$ is a set of isolated points and $p\in M^T$,
then $A_{0\dots 0}(p)=1$ and the others $A_{n_1\dots n_s}(p)$
vanish. Since $C_-(0,\,l)=1$, $\,C_+(0,\,l)=1$ for $l>0$ and
$C_+(0,\,0)=0$ the number of nonzero summands in
 $$\sum_{l\in {\mathcal L}(p,k)}
 \big(\prod_{j=1}^n \tilde C_{\sigma_{pj}}(0,l_j)\big)$$
 is precisely the number $N_p(k,u)$ introduced in Subsection
 \ref{SectionT}. Thus we deduce from (\ref{Formula})
 $$\sum_{p\in Q_+(u)}N_p(k,u)-\sum_{p\in
 Q_-(u)}N_p(k,u)=0,$$
which agrees with  Theorem
 \ref{ThmGeneralGeneral}.

\smallskip

{\it Example 2.}

 Let us suppose that ${\rm dim}\,M=4$ and  it is acted on by $S^1$ so that
 $M^{S^1}$ has two connected components: A point $q$ and  $F$, a
 $2$-submanifold of $M$.
Let us assume that $J(q)$ is the    minimum  value of $J$ and
$J(F)$ the maximum one. So $\alpha_{qj}>0$ for $j=1,2$ and
$\alpha_F<0$.

It is easy to determine the contribution of $q$ to (\ref{Formula})
since the only nonzero $A_{n_1n_2}(q)$ is $A_{00}(q)=1$. This
contribution is
$$\sharp\{(l_1,\,l_2)\,|\, l_j\in {\mathbb N}_{>0}, \,
\alpha_{q1}l_1+\alpha_{q2}l_2=k-J(q)\}.$$

On the other hand $A_0(F)=\int_F(\omega+{\rm Td}_1(F))$ and
$A_1(F)=-\int_Fc_1({\mathcal N}_F)$. The other $A_n(F)$ are equal
to zero. Hence
$$\chi(t)=
\sum_{m_1,m_2\geq 1}t^{m_1\alpha_{q1} + m_2\alpha_{q2}+J(q)}
+A_0\sum_{l\geq 0}t^{J(F)-l\alpha_F} +A_1 \sum_{l_1,l_2\geq 0}
t^{J(F)-(l_1+l_2+1)\alpha_F}.
$$
For any $m_1>>1$ the exponent of the  monomial
$t^{m_1\alpha_{q1}+\alpha_{q_2}+J(q)}$
must appear in a monomial of type $t^{J(F)-m\alpha_F}$. Therefore
$J(F)-J(q)\in I_-$, where $I_-$ is the ideal of ${\mathbb Z}$
generated by $\alpha_F$.

Given $k\in{\mathbb N},$ if
 \begin{equation}\label{n0equation}
 n_0:=\frac{J(F)-k}{\alpha_F}\notin{\mathbb N},
 \end{equation}
the coefficient of $t^k$ in $\chi(t)$ is
$$N_q(k)=\sharp\{(m_1,\,m_2)\,|\, m_j>0,\;m_1\alpha_{q1}+m_2\alpha_{q2}+J(q)=k \}=0.$$
Thus for $k$ sufficiently large, $N_q(k)=0$ (assumed that
(\ref{n0equation}) holds).

If $ n_0:=(J(F)-k)/{\alpha_F}\in{\mathbb N}$, the coefficient of
$t^k$ in $\chi(t)$ is $N_q(k)+A_0+n_0A_1$, since
$$\sharp\{(l_1,\,l_2)\,|\,l_j\in{\mathbb N},\; J(F)-(l_1+l_2+1)\alpha_F=k  \}=n_0.$$
So we can state the following Proposition.

%
 \begin{Prop}\label{sharpformula} Let $M$ be a $4$-manifold acted on by
 $S^1$.
If   the fixed point set has only two
 components, $q$ and $F$ of dimensions $0$ and $2$ respectively and  $J(F)>J(q)$, then

1) There is a natural number $m$ such mapping $$ n_0\in{\mathbb
N}_{>m}\mapsto\sharp\{(m_1,\,m_2)\,|\, m_j>0,\; m_1\alpha_{q1}+
m_2\alpha_{q2}=J(F)-J(q)-n_0\alpha_F \}$$
 is the affine map $-A_0(F)-n_0A_1(F).$

2)
$$\sharp\{(m_1,\,m_2)\,|\, m_j>0,\;m_1\alpha_{q1}+m_2\alpha_{q2}=k-J(q) \}=0$$
for all $k$ sufficiently large such that
$(J(F)-k)/\alpha_F\notin{\mathbb N}$.

3)\; $J(F)-J(q)\in I_-,$ where
 $I_-$ is the ideal of ${\mathbb Z}$ generated by $\alpha_F$.
 \end{Prop}

\smallskip

Now we return to the toric manifold $M_{\Delta}$ considered in
{\it Example 1} (Subsection (\ref{SectionU(1)})).
$$M_{\Delta}=\{(z_0,z_1,z_2)\,|\, \sum_j|z|^2=1  \}/\sim.$$
 The $T=(U(1))^2$-action on $M_{\Delta}$ is defined by (see
 \cite{Gui})
 $$(\lambda_0,\,\lambda_1)[z]=[\lambda_0z_0,\lambda_1z_1,z_2].$$
 We write $z_j=\rho_j\,{\rm exp}(i\theta_j)$, with $\theta_j\in{\mathbb
 R}/{\mathbb Z}$. Then $(\rho_0,\varphi_0,\rho_1,\varphi_1)$, with
 $\varphi_j=\theta_j-\theta_2$, are coordinates on $M_{\Delta}$,
 and the map $J_T([z])=(\rho_0^2,\,\rho_1^2)$ satisfies ${\rm
 im}\,J_T=\Delta$. If $Y=(ai, bi)\in{\frak t}$, then
 $Y_M=-a\frac{\partial}{\partial\varphi_0}-b\frac{\partial}{\partial\varphi_1}.$
 Our convention $\omega(Y_M,\,)=d\langle\mu,\,Y\rangle=i
 d\langle  J_T,\,Y\rangle$ imposes that
 $$\omega=-\big(d\rho_0^2\wedge d\varphi_0+d\rho_1^2\wedge d\varphi_1  \big).$$

 Next we consider the $S^1$-action on $M_{\Delta}$ generated by
 $X=(x,\,x)\in{\mathbb Z}^2$. The components of the fixed point set
 for this $S^1$-action are $p=[0,0,1]$ and $F=\{[z]\,|\, z_2=0\}\simeq {\mathbb C}P^1.$
 The weights of the isotropy representations are
 $$\alpha_{p1}=x,\; \alpha_{p2}=x,\; \alpha_{F}=-x,\;\,\text{and }\, J(p)=0,\;J(F)=x.$$
 The tangent bundle $TF$ is ${\mathcal O}(2)$ and the normal bundle of
 ${\mathbb C}P^1$ in ${\mathbb C}P^2$ is ${\mathcal O}(1)$.
 Moreover
 $$\int_F\omega=-\int_{0}^1d\,\rho_0^2\int_0^1d\,\varphi_0=-1.$$
 Thus $A_0=0$, $A_1=-1$ and $-A_0(F)-n_0A_1(F)=n_0$. Furthermore
  $\sharp\{(m_1,m_2)\,|\,m_j>0,\; m_1+m_2=1+n_0
 \}=n_0$, in accordance with the first item of Proposition
 \ref{sharpformula}.

\smallskip



\begin{thebibliography}{99}








\bibitem{A}
 Atiyah, M.F.: Convexity and commuting Hamiltonians. Bull. London
 Math. Soc. {\bf 14}, 1-15 (1982).




\bibitem{A-S}
  Atiyah, M. F.,  Segal G.B¨: The index of elliptic operators II. Ann. of Math {\bf 87}, 531-545 (1968).






\bibitem{B-G-V}
  Berline, N.,  Getzler, E.,  Vergne,  M.:
   {\it Heat kernels and Dirac operators}.
Springer-Verlag, Berlin, 1991.


\bibitem{Brion-V}
Brion, M.,  Vergne, M.: Residue formulae, vector partition
functions and lattice points in rational polytopes.   J. Amer.
Math. Soc. {\bf 10}, 797-833 (1997).




\bibitem{D-G-M-W}
Duistermaat, H., Guillemin, V., Meinrenken, E.,  Wu, S. :
Symplectic reduction and Riemann-Roch for circle actions.
 Math. Res. Lett. {\bf 2}, 259-266  (1995).









\bibitem{G}
 Guillemin,   V.:
 Reduced phase spaces and
Riemann-Roch. Lie theory and geometry, 305-334, Progr. Math., 123,
Birkhäuser Boston, Boston, MA, 1994.


\bibitem{Gui}
 Guillemin,   V.:
  {\it Moment maps and combinatorial invariants of Hamiltonian
   $T^n$-spaces}.
  Birkh\"auser,
  Boston,
  1994.



\bibitem{Gui1}
   Guillemin, V.: Riemann-Roch for toric orbifolds
     J. Differential Geom. {\bf 45}, 53-73 (1997).





\bibitem{G-S-0}
   Guillemin, V.,   Sternberg, S.:
   Convexity properties of the moment mapping.   Invent. Math. {\bf 67}, 491-513 (1982).



\bibitem{G-S}
   Guillemin, V.,   Sternberg, S.:
   Geometric quantization and multiplicities of group
   representations. Invent. Math. {\bf 67}, 515-538 (1982).





\bibitem{J-K}
Jeffrey, L.C., Kirwan, F.C.:  On localization and Riemann-Roch
numbers for symplectic quotients.
 Quart. J. Math. Oxford Ser. (2)
{\bf 47}, no. 186, 165--185 (1996).


\bibitem{K}
 Kirwan, F.: Convexity properties of the moment mapping, III.
  Invent. math. {\bf 77}, 547-552 (1984).



\bibitem{L-M}
 Lawson, H. B., Michelsohn, M.-L. :
 {\it Spin geometry}.  Princeton University Press,
 Princeton, 1989.



\bibitem{Mc-S}
 McDuff,  D.,   Salamon, D.:
 {\it Introduction to symplectic topology}. Clarenton Press,
 Oxford,
1998.



\bibitem{Mein}
Meinrenken, E.: Symplectic surgery and $\text{Spin}^c$-Dirac
operator. Adv. Math. {\bf 134}, 240-277 (1998).




\bibitem{lP01}
 Polterovich,   L.:{\it The geometry of the group of symplectic
diffeomorphisms}. Birkh\"auser, Basel, 2001.

\bibitem{SJ}
   Sjamar,  R.:
   Symplectic reduction and Riemann-Roch formulas for
   multiplicities.
    Bull. Amer. Math. Soc. (N.S.)
  {\bf 33}, 327-338, (1996).

  \bibitem{Verg}
 Vergne, M.: Multiplicities formulas for geometric quatization.
 Part I, II. Duke Math. J. {\bf 82} 143-194 (1996).

\bibitem{V}
   Vi\~na,  A.:
   Symplectic action around loops in ${\rm Ham}(M).$
   Geom. Dedicata {\bf 109}, 31-49, (2004).



\bibitem{We}
 Weinstein,  A.:
   Cohomology of symplectomorphism groups and critical values of
   Hamiltonians.
    Math. Z. {\bf 201}, 75-82 (1989)



\bibitem{W}
 Woodhouse, N.M.J.  {\it Geometric quantization}. Clarenton Press, Oxford, 1992.















\end{thebibliography}
\end{document}